\documentclass[a4paper, 10pt]{amsart}
\usepackage{amsfonts,mathrsfs,amsmath,amscd,amssymb}
\usepackage[pdftex]{graphicx}
\usepackage{hyperref}
\usepackage{color}
\newtheorem{theorem}{Theorem}[section]
\usepackage[utf8]{inputenc}
\usepackage{enumitem}
\usepackage{mathtools}
\usepackage[english]{babel}
\usepackage[all,cmtip]{xy}
\usepackage{soul} 

\newtheorem{proposition}[theorem]{\sf Proposition}
\newtheorem{lemma}[theorem]{\sf Lemma}
\newtheorem{definition}[theorem]{\sf Definition}
\newtheorem{corollary}[theorem]{\sf Corollary}
\newtheorem{remark}[theorem]{\sf Remark}

\def\C{\mathbb C}
\def \Z{\mathbb Z}
\def \Q{\mathbb Q}
\def \R{\mathbb R}

\def \O{\mathcal O}
\def \A{\mathcal A}
\def \X{\mathfrak{X}}
\def \E{\mathcal{E}}
\def \L{\mathcal{L}}
\def \F{\mathcal{F}}
\def \G{\mathbb G}
\def \mcG{\mathcal G}
\def \ra{\rightarrow}
\def \mcX{\mathcal{X}}
\def \Y{\mathcal{Y}}

\let\olddefinition\definition
\renewcommand{\definition}{\olddefinition\normalfont}

\let\oldremark\remark
\renewcommand{\remark}{\oldremark\normalfont}

\DeclareMathOperator{\Br}{Br}

\DeclareMathOperator{\Spec}{Spec}

\DeclareMathOperator{\Supp}{Supp}

\begin{document}
\author{Tanya Kaushal Srivastava}
\address{Indian Institute of Technology, Gandhinagar, Gandhinagar-382355, India}.
\email{ks.tanya@iitgn.ac.in}
\author{Sofia Tirabassi}
\address{Stockholm University, 106 91 Stockholm, Sweden}
\email{tirabassi@math.su.se}
\title[Tame Twisted FM partners]{Counting  Twisted Tame Fourier-Mukai Partners of an Ordinary K3 Surface}

\maketitle
\tableofcontents

\begin{abstract}
    In this article, we prove that a tame twisted K3 surface over a algebraically closed field of positive characteristic has only finitely many tame twisted Fourier-Mukai partners and we also give a counting formula in case we have an ordinary tame untwisted K3 surface.
\end{abstract}
\section{Introduction}

The goal of this article is to prove that the twisted tame Fourier-Mukai partners of a K3 surface over an algebraically closed field of positive characteristic are finitely many and to count the number of twisted tame Fourier-Mukai partners of an ordinary K3 surface. More precisely, we show that 
\begin{theorem}
Let $X$ be a K3 surface over $k$, an algebraically closed field of positive characteristic greater than 3, and $\alpha \in Br(X)$, then the number of pairs $(Y, \beta)$, where $Y$ is a K3 surface and $\beta \in Br(Y)$ is a Brauer class such that $D^b(X, \alpha) \cong D^b(Y, \beta)$, with orders of $\alpha$ and $\beta$ are not divisible by the characteristic of $k$, is finite and every such pair $(Y, \beta)$ comes from a coarse moduli space of stable twisted sheaves on $X$ with natural twisting. 
\end{theorem}

In case we start with an untwisted ordinary K3 surface, we can even count the number of tame Fourier-Mukai partners as follows: 
\begin{theorem}
Let $X$ be a ordinary K3 surface. Then the number of twisted tame Fourier-Mukai partners of $X$ is the same as the number of twisted Fourier-Mukai partners of $X_{can, \bar{K}}$, the geometric generic fiber of the canonical lift. 
\end{theorem}

The number of twisted Fourier Mukai partners of an untwisted K3 surface defined over field of complex numbers is given by the following result of Ma.

\begin{theorem}[\cite{MaFM} Theorem 4.2]
Let $S$ be a projective K3 surface over $\C$ with transcendental lattice $T(S)$, Neron-Severi lattice $NS(S)$ of discriminant $D_{NS(S)}$ and $FM^d(S)$ the set of isomorphism classes of twisted Fourier-Mukai partners $(S', \alpha')$ of $S$, where $\alpha'$ is an element of order $d$ in the Brauer group of S. Then the number of such isomorphism classes is given by the following formula:

\begin{eqnarray*}
\# FM^d(S) = &\sum_x \Big\{ \sum_M \#\Big(O_{Hodge}(T_x, \alpha_x) \backslash O(D_M)/O(M)\Big) +  \\
& \epsilon(d)\sum_{M'}\#\Big( O_{Hodge}(T_x, \alpha_x) \backslash O(D_{M'})/O(M')  \Big) \Big\},
\end{eqnarray*}
where $x$ runs over the set $O_{Hodge}(T(S))\backslash I^d(D_{NS(S)})$ ($I^d(D_{NS(S)})$ is the set of order $d$ isotropic elements of the discriminant form $D_{NS(S)}$), the lattices $M, M'$ run over the certain subsets ($\mathcal{G}_1, \mathcal{G}_2$, depending on $x$) of the genus set $\mathcal{G}$, and $\epsilon(d)$ is defined to be $1$ if $d=1,2$ and $2$ if $d\geq 3$.  

\end{theorem}

Note that the above counting formula is for a given order of the Brauer class, so we can split the case of counting twisted Fourier-Mukai partners of ordinary K3 surfaces in two broad classes: one in which the order is not divisible by the characteristic, tame classes, and the other in which it is, wild classes. In this article we will tackle the tame classes.  The wild cases will be discussed in our next article. 

\subsubsection*{Plan of the Article}

The proof of our main theorem relies on our ability to lift a K3 surface over a field of characteristic $p$ to characteristic $0$ along with its Brauer class and having a special lift (canonical lift of an ordinary K3 surface) which we know how to identify from among all the lifts of a K3 surface. 

We begin in section 2 by recalling the definition and basic results about the Brauer group of a K3 surface. We then study the lifting of a (tame) Brauer class from characteristic $p$ to $0$ and prove that the lifts are unobstructed and unique. These results rely on tameness and to extend our results to tame case, one needs to develop the arguments for the wild case.  

In section 3, we recall the definition of an ordinary K3 surface using the notion of height of a F-isocrsytal and then define their canonical lift and state the Taelman criterion to identify canonical lifts. 

In section 4 lies the major part of the proof of our main theorem. This is where we prove that every twisted Fourier-Mukai partner of a twisted K3 surface is a moduli space of twisted sheaves (Theorem \ref{TwistedpartnersModuli}). We begin by recalling the definitions of a Twisted K3 surface and twisted Fourier-Mukai partner, the construction of moduli space of twisted sheaves. To be able to prove Theorem \ref{TwistedpartnersModuli}, we needed to extend the lifting result of Lieblich-Olsson \cite[Theorem 6.3]{LO}  to the twisted K3 surfaces and this result is proved as Theorem \ref{LOlifting}. 

In the last section, we end by proving the finiteness and counting results. 

\subsection*{Notations and Conventions}

Unless stated otherwise we will follow the following convention:
\begin{table}[hbt!]
\begin{tabular}{l l}
$k$ & Algebraically closed field of positive characteristic \\
$W(k)$ &  Ring of Witt vectors of $k$ \\
 $Frob_W$ & Frobenius morphism of $W(k)$\\
 $K$ & Field of fractions of $W(k)$ \\
 $Frob_K$ & Frobenius morphism of $K$.  
\end{tabular}
\end{table}

\section{Brauer group and deformations}

In this section, we start by recollecting all the results about the Brauer group of a K3 surface that we will be using in the further sections mostly rather implicitly. Here we will also prove the existence of a unique lift of a tame Brauer class on a K3 surface from characteristic $p$ to characteristic $0$. 

\subsection{Brauer Group of a K3 surface}

In this subsection, we recall various definitions of Brauer groups associated to a scheme and the connections between them for a K3 surface. We refer the reader to \cite[Chapter 3]{CTS} and \cite[Chapter 18]{HuyK3} for more details.

\subsubsection{Brauer-Azumaya group}
\begin{definition}[\cite{Gro68b}, Theorem 5.1]
An \textbf{Azumaya algebra} on a scheme $X$ is an $\O_X$-algebra $\mathcal{A}$ that is coherent as an $\O_X$-module with $\A_x \neq 0$ for all $x \in X$, and that satisfies one of the following equivalent conditions:
\begin{enumerate}
\item There is an open covering $\{U_i \ra X \}$ in the \'etale topology such that for each $i$ there exists $r_i \in \Z_{\geq 0}$ such that $\A \otimes_{O_X} \O_{U_i} \cong M_{r_i}(\O_{U_i})$
\item $\A$ is locally free as an $\O_X$-module and the canonical homomorphism $\A \otimes_{\O_X} \A^{op} \ra \mathcal{E}nd_{\O_X}(\A)$ is an isomorphism.
\end{enumerate} 

The minimal of the $r_i$ in the above point (1) is called the index of the Azumaya algebra.  
\end{definition}

\begin{definition}
An Azumaya algebra $\A$ on $X$ is called \textbf{trivial} if it is isomorphic to endomorphism algebra of a locally free sheaf, i.e., $\A \cong \mathcal{E}nd_X(E)$, where $E$ is a vector bundle on $X$. Two Azumaya algebras $\A$ and $\A'$ are called \textbf{equivalent} if there exist vector bundles $E$ and $E'$ such that $\A \otimes_A \mathcal{E}nd(E) \cong \A' \otimes \mathcal{E}nd(E')$.   
\end{definition}

\begin{definition}
The \textbf{Brauer group} of a scheme $X$ is the group of Azumaya algebra on $X$ upto equivalence. The inverse of an Azumaya algebra $A$ is given by the opposite Azumaya algebra $\A^{op}$ as $\A \otimes \A^{op} \cong \mathcal{E}nd{A}$.    \end{definition}

We will denote this group as $Br(X)$. 

Note that set of Azumaya algebras of index $n$ upto equivalence is in 1-1 correspondence with $H^1_{et}(X, PGL_n)$.   

\subsubsection{Cohomological Brauer group}

Via the following short exact sequence 
$$
0 \ra \G_m \ra GL_n \ra PGL_n \ra 0
$$
we see that there is there is a natural group homomorphism (via the boundary homomorphism of the corresponding cohomological long exact sequence).
$$
Br(X) \ra H^2(X, \G_m).
$$
The second group was used by Grothendieck to define:  

\begin{definition}
The \textbf{cohomological Brauer group} for a scheme $X$, denoted as $Br'(X)$ is defined to be:
$$
Br'(X) = H^2_{et}(X, \G_m)_{tors}.
$$
\end{definition}

In case the underlying scheme is quasi-compact and separated with an ample line bundle then it is a result of Gabber \cite[4.2.1]{CTS} that the two Brauer groups define above are isomorphic, i.e., 
$$
Br(X) \cong \Br'(X). 
$$

In case the underlying scheme is a regular, integral, noetherian scheme we have that
$$
H^2_{et}(X, \G_m) = H^2_{et}(X, \G_m)_{tors},
$$
i.e., the second etale cohomology group of $X$ with coefficients in $G_m$ is a torsion group. See \cite[Lemma 3.5.3]{CTS} for a more general statement. 

Thus, for a K3 surface we have the following isomorphism:
$$
H^2(X, \G_m) \cong Br'(X) \cong Br(X).
$$

\subsubsection{Formal Brauer group}
As one defines the Picard functor for a scheme, we can define a Brauer group functor, it unfortunately turns out that this functor is not representable by an algebraic group. However, its formal completion turns out to be a formal group scheme, as constructed by Artin and Mazur, see \cite[Section 18.1.3]{HuyK3} for more details and references, we give the definition here so that we can compare with the definition below and for an equivalent definition of an ordinary K3 surface.

Let $(Art/k)$ denote the category of Artin local k-algebras, $(Ab)$ the category of Abelian groups and $X_A$ the base change of $X$ to $\Spec (A)$, for $A \in (Art/k)$
\begin{definition}
Consider the functor:
\begin{eqnarray*}
\hat{\textbf{Br}_X}: (Art/k) &\ra &(Ab) \\
A &\mapsto  & Ker(Br(X_A) \ra Br(X)).
\end{eqnarray*}
If $\hat{\textbf{Br}_X}$ is prorepresentable by a complete $k$-algebra $R$, we call the formal group $Spf(R)$, the \textbf{formal Brauer group} of $X$, denoted by $\hat{Br_X}$.
\end{definition}

For a K3 surface over any field $k$, the functor $\hat{\textbf{Br}_X}$ is pro-representable by a smooth, one dimensional formal group. 

\subsubsection{Enlarged Formal Brauer group} \label{enlargeFBG}

Artin and Mazur \cite{AM} enlarged the above formal Brauer group functor by considering the fppf group cohomology with coefficients in $\mu_{p^{\infty}}$. To be precise, let $A$ be an Artin local ring with residue field $k$, $\mathcal{X}/ A$ be a lifting of the scheme $X/ k$ and let $(Art/ A)$ be the category of Artin local A-algebras with residue field $k$, then we can define the following:
\begin{definition}
The \textbf{enlarged formal group functor} is defined as 
\begin{eqnarray*}
\Psi_{\mathcal{X}}: (Art/A) &\ra &(Ab) \\
B &\mapsto & H^2_{fppf}(\mathcal{X}_B, \mu_{p^{\infty}})
\end{eqnarray*}
If this functor is representable by a group scheme, we call it the \textbf{enlarged formal Brauer group} of $\mcX$.
\end{definition}

Compare with the remarks in \cite[Section 2.1]{Taelman}. These groups will be very useful to us in giving the definition of canonical lift of an ordinary K3 surface in Section \ref{Section 3}. In the situation of ordinary K3 surfaces (defined below), we note that the formal Brauer group functor is representable by a p-divisible group, see \cite[Theorem 1.6]{Taelman}.

\subsection{Deformation of Azumaya Algebra} \label{2.2}
(This argument uses tameness.)
Let $X$ be a K3 surface over an algebraically closed field $k$ of positive characteristic. The obstruction to deforming an Azumaya algebra $\A$ lies in the second cohomology group $H^2(X,  I \otimes (\A/\O_X))$, where $I$ is the sheaf of ideals of a square zero extension $X \hookrightarrow X'$ \cite[Lemma 3.1]{deJong}. When unobstructed, the deformations of the Azumaya algebra are in one-to-one to correspondence with the set $H^1(X, I \otimes (\A/\O_X))$. 
In case characteristic of $k$ does not divide $n^2$, the rank of $\A$, de Jong \cite[Proposition 3.2]{deJong} has shown that can always modify the Azumaya algebra to become unobstructed. More precisely, there exists an elementary transform $A'$ of A such that they have the same classes in the Brauer group of $X$ and $H^0(X, K\A) = 0$. 
The above discussion shows that given a Brauer class on a K3 surface (the order of Brauer class should not be divisible by the characteristic of $k$) one can always pick an unobstructed Azumaya algebra representing it, thereby lifting the Brauer class itself to characteristic zero.

\subsection{Deformation of Brauer class on K3 surfaces} \label{2.4}

Let $X$ be a K3 surface over a field of positive characteristic $p$. Let $H^2(X, \O_X^*)$ be the Brauer group of $X$, for any small extension 
$$
0 \ra \mathfrak{m} \ra A' \ra A \ra 0
$$
of local Artin $W$-algebras, we have the following short exact sequence:
$$
0 \ra \mathfrak{m}\otimes_k \O_X \ra \O_{X_{A'}} \ra \O_{X_A} \ra 0, 
$$
where $X_{A'}$ and $X_{A}$ are the corresponding lifts of X. 
Restricting to invertible elements we get
$$
0 \ra 1+\mathfrak{m}\otimes_k\O_X \ra \O_{X_{A'}}^* \ra \O_{X_{A}}^* \ra 0
$$
This induces the following long exact sequence of cohomology:
\begin{eqnarray*}
\ldots \ra H^1(X, 1 + \mathfrak{m}\otimes_k \O_X) \ra H^1(X_{A'},\O_{X_{A'}}^*) \ra H^1(X_A, \O_{X_{A}}^*) \ra \\
H^2(X, 1 + \mathfrak{m}\otimes_k \O_X) \ra H^2(X_{A'},\O_{X_{A'}}^*) \ra H^2(X_A, \O_{X_{A}}^*) \ra \\
H^3(X, 1 + \mathfrak{m}\otimes_k \O_X) \ra \ldots. 
\end{eqnarray*}

Then 
$$H^3(X, 1 + \mathfrak{m}\otimes_k \O_X) = H^3(X, \O_X) 
\otimes_k (1 + \mathfrak{m}) = 0 
$$ as for any surface we have $H^3(X, \O_X) = 0$. 
Thus, deformation of Brauer class is unobstructed and we can lift them. 

\subsubsection{Algebraization of Lifted Brauer classes}

Using the results of the previous section, we can now construct a formal lift of our K3 surface along with its Brauer class to characteristic zero. Moreover, this formal lift is algebraizable. 

Recall that Grothendieck Existence theorem \cite[Theorem 21.2]{HartDef} implies that the Picard preserving\footnote{Actually, it is enough that one ample line bundle lifts, but since we will always be working with lifts that preserve the Picard group, we use this statement. Canonical lifts are example of Picard preserving lifts.} formal lifts of K3 surfaces are algebraic.Thus we have a smooth projective scheme $\mathcal{X}$ over the Witt ring $W(k)$. 

As in our case the order of the Brauer class $\alpha$ is coprime to characteristic of the base field $p$, then we can use deJong's argument above to choose an unobstructed Azumaya algebra associated to $\alpha$ and then the algebraization of the lifted Brauer class is the same as that of lifted Azumaya algebra and the algebraization follows from Grothendieck existence theorem.

\subsubsection{Uniqueness of lifted Brauer class}

Again using the long exact sequence in section \ref{2.4} above, we note that the lifts will be unique if the maps 
$$
H^2(X, 1 + \mathfrak{m}\otimes_k \O_X) \ra H^2(X_{A'},\O_{X_{A'}}^*)
$$
are zero. 
Now since we in tame case, note that all we have to show is that there is no prime to p torsion elements in $H^2(X, \O_X)$. Recall that for a K3 surface this is just $k$ as a group and hence has only p-torsion and we are done.

\section{Ordinary K3 surfaces} \label{Section 3}

Over fields of positive surfaces, K3 surfaces can be broadly put in three classes, depending on the notion of height of a K3 surface, namely as:
\begin{enumerate}
    \item Height 0: Supersingular K3 surfaces. By convention the height is sometimes also said to be infinite.
    \item Height 1: Ordinary K3 surfaces
    \item Height greater than 1 but finite: K3 surfaces of finite height.
\end{enumerate}

Out of these classes, the ordinary K3 surfaces behave mostly like the K3 surfaces of characteristic, essentially due to the fact that the admit a canonical lift to characteristic zero. Moreover, there exist cohomological criterion to identify canonical lift of an ordinary K3 from its other lifts. 

On the other hand supersingular K3 surfaces display most of the pathological positive characteristic behaviour. 

Finite height K3 surfaces are expected to behave more like K3 surfaces over fields of characteristic zero but they don't admit canonical lift thereby making direct comparison methods to study them infeasible. 

In this section, we begin by by recalling the notion of height and then we recall some results about ordinary K3 surfaces and their canonical lifts as proved by Nygaard in \cite{Nygaardtate} and \cite{Nygaardtorelli}, and by Deligne-Illusie in \cite{Deligne-Illusie}. We will define the height of a K3 surface through its F-crystal, since this is the characterisation we will be using later in section \ref{HeightTDI} to prove that height is a twisted derived invariant. 

For the definition of height via the formal Brauer groups see \cite{HuyK3} and \cite{Liedtke}. Both definitions turn out to be equivalent (for example see Prop. 6.17 \cite{Liedtke}).

\begin{definition}[F-(iso)crystal] \label{F-crystal}
An \textbf{F-crystal} $(M, \phi_M)$ over $k$ is a free $W$-module $M$ of finite rank together with an injective $Frob_W$-linear map $\phi_M: M \ra M$, that is, $\phi_M$ is additive, injective and satisfies 
$$
\phi_M (r \cdot m) = Frob_W(r) \cdot \phi_M(m) \ \text{for all} \ r \in W(k), m \in M.
$$
An \textbf{F-isocrystal} $(V, \phi_V)$ is a finite dimensional $K$-vector space $V$ together with an injective $Frob_K$-linear map $\phi_V: V \ra V$. 

A \textbf{morphism $u:(M, \phi_M) \ra (N, \phi_N)$ of F-crystals} (resp. \textbf{F-isocrystals}) is a $W(k)$-linear (resp. $K$-linear) map $M \ra N$ such that $\phi_N \circ u = u \circ \phi_M$. An \textbf{isogeny} of F-crystals is a morphism $u: (M, \phi_M) \ra (N, \phi_N)$ of F-crystals, such that the induced map $ u \otimes Id_K: M \otimes_{W(k)} K \ra N \otimes_{W(k)} K$ is an isomorphism of F-isocrystals.
\end{definition}

\textbf{Examples:}
\begin{enumerate}
\item The trivial $F$-crystal: $(W, Frob_W)$.
\item  Geometric $F$-crystal: Let $X$ be a smooth and proper variety over $k$. For any $n$, take the free $W(k)$ module $M$ to be $H^n := H^n_{crys}(X/W(k))/torsion$ and $\phi_M$ to be the Frobenius $F^*$. The Poincar\'e duality induces a perfect pairing 
\begin{equation*}
\langle -, - \rangle: H^n \times H^{2dim(X)-n} \ra H^{2dim(X)} \cong W 
\end{equation*}
which satisfies the following compatibility with Frobenius 
\begin{equation*}
\langle F^*(x), F^*(y) \rangle = p^{dim(X)} Frob_W(\langle x, y \rangle),
\end{equation*}
where $x \in H^n$ and $y \in H^{2dim(X)-n}$. As $Frob_W$ is injective, we have that $F^*$ is injective. Thus, $(H^n, F^*)$ is an F-crystal. We will denote the F-isocrystal $H^n_{crys}(X/W) \otimes K$ by $H^n_{crys}(X/K)$.
\item The F-isocrystal $K(1) : = (K, Frob_K/p)$. Or more generally, one has the F-isocrystal $K(n) := (K, Frob_K/p^n)$ for all $n \in \Z$. Moreover, for any F-isocrystal $V$ and $n \in \Z$, we denote by $V(n)$ the F-isocrystal $V \otimes K(n)$. 
\end{enumerate}

 Recall that a theorem of Dieudonn\'e-Manin gives us that the category of F-crystals over $k$ up to isogeny is semi-simple and the simple objects are the F-crystals: 
\begin{equation*}
M_{\alpha} = ((\Z_p[T])/(T^s-p^r)) \otimes_{Z_p} W(k), (\text{mult. by} \ T)\otimes Frob_W),
\end{equation*} 
for $\alpha = r/s \in \Q_{\geq 0}$ and $r$, $s$ non-negative coprime integers. 

We define the \textbf{rank} of the F-crystal $M_{\alpha}$ as $s$. We call $\alpha$ the \textbf{slope} of the F-crystal $M_{\alpha}$. 

\begin{definition} Let $(M, \phi)$ be an F-crystal over $k$ and let 
$$
(M, \phi) \sim^{isogeny} \oplus_{\alpha \in \Q_{\geq 0}} M_{\alpha}^{n_{\alpha}}
$$
be its decomposition up to isogeny. Then the elements of the set 
$$
\{\alpha \in \Q_{\geq 0}| n_{\alpha} \neq 0 \}
$$
are called the \textbf{slopes} of $(M, \phi)$. For every slope $\alpha$ of $(M, \phi)$, the integer $\lambda_{\alpha} := n_{\alpha} \cdot rank_W M_{\alpha}$ is called the \textbf{multiplicity} of the slope $\alpha$.
\end{definition}

Moreover, the above classification result of Dieudonn\'e-Manin also gives that any F-isocrystal $V$ with bijective $\phi_V$ is isomorphic to a direct sum of F-isocrystals
$$
(V_{\alpha} := K[T]/(T^s-p^r), (\text{mult. by $T$}) \otimes Frob_K), 
$$   
for $\alpha = r/s \in \Q$. The dimension of $V_{\alpha}$ is $s$ and we call $\alpha$ the \textbf{slope} of $V_{\alpha}$.

\begin{definition}[Height] \label{heightdef}
The \textbf{height} of a K3 surface $X$ over $k$ is the sum of multiplicities of slope strictly less than 1 part of the F-crystal $H^2_{crys}(X/W)$. In other words, the dimension of the subspace of slope strictly less than one of the F-isocrystal $H^2_{crys}(X/K)$, which is $\text{dim} (H^2_{crys}(X/K)_{[0,1)} := \oplus_{\alpha_i < 1}V_{\alpha_i}^{n_{\alpha_i}})$.
\end{definition}

If for a K3 surface $X$ the $\text{dim} (H^2_{crys}(X/K)_{[0,1)}) =0$, then we say that the height of $X$ is infinite. A K3 surface with infinite height is called a \textbf{supersingular K3 surface}.

\begin{definition}[Ordinary K3 surface]
A K3 surface $X$ over a perfect field $k$ of positive characteristic is called \textbf{ordinary} if the height of $X$ is 1. 
\end{definition}

\begin{proposition}[\cite{Nygaardtate} Lemma 1.3] The following are equivalent:
\begin{enumerate}
\item $X$ is an ordinary K3 surface,
\item The height of formal Brauer group is 1,
\end{enumerate}
\end{proposition}

Thus, one can even use formal Brauer group to define an ordinary K3 surface. We now give a quick overview about canonical lift of an ordinary K3 surface.

Recall the definition of enlarged formal group from section \ref{enlargeFBG} for an ordinary K3 surface. In \cite{AM}  Artin-Mazur showed that the enlarged Brauer group $\Psi_{X_A}$ defines a $p$-divisible group on $\text{Spec} (A)$ lifting $\Psi_X / k$. 

\begin{theorem}[Nygaard \cite{Nygaardtate}, Theorem 1.3]
Let $X/k$ be an ordinary K3 surface. The map
\begin{equation*}
\{ \text{Iso. classes of liftings $X_A/A$} \} \ra \{ \text{Iso. classes of liftings G/A} \}
\end{equation*}
defined by 
\begin{equation*}
X_A/A \mapsto \Psi_{X_A}/A
\end{equation*}
is a functorial isomorphism. 
\end{theorem}

From \cite[Proposition IV.1.8]{AM}, we know that  the enlarged Brauer group of an ordinary K3 surface fits in the following exact sequence
\begin{equation}
0 \ra \Psi_X^0 (= \hat{Br_X}) \ra \Psi_X \ra \Psi^{\acute{e}t} \ra 0.
\end{equation}
 
As the height one formal groups are rigid, there is a unique lifting $G^0_A$ of $\Psi_X^0$ to $A$. Similarly, the \'etale groups are rigid as well, hence there is a unique lift $G^{\acute{e}t}_A$ of $\Psi_X^{\acute{e}t}$ to $A$. This implies that if $G$ is any lifting of $\Psi_X$ to $A$, then we have an extension
$$
0 \ra G^0_A \ra G \ra G^{\acute{e}t}_A \ra 0
$$
lifting the extension
$$
0 \ra \Psi_X^0 \ra \Psi_X \ra \Psi_X^{\acute{e}t} \ra 0.
$$
Therefore, the trivial extension $G = G^0_A \times G_A^{\acute{e}t}$ defines a unique lift $X_{can,A}/A$  of $X/k$ such that $\Psi_{X_{can,A}} = G^0_A \times G^{\acute{e}t}_A$. Take $A = W_n$ and $X_n = X_{can, W_n}$, thus we get a proper flat formal scheme $\{X_n\}/ Spf W$.

\begin{theorem}[Definition of \textbf{Canonical Lift}]
The formal scheme $\{X_n\} / Spf W$ is algebraizable and defines a K3 surface $X_{can}/ \Spec(W)$. 
\end{theorem}
This is a result of Nygaard \cite[Proposition 1.6]{Nygaardtate}.

One of the most important property of the canonical lift is that it is a Picard lattice preserving lift, i.e., The canonical lift $X_{can}$ has the property that any line bundle on $X$ lifts uniquely to $X_{can}$, \cite[Proposition 1.8]{Nygaardtate}.

Lastly, we state a criteria to determine when a  lifting of an ordinary K3 surface is the canonical lift. This is the criteria that we will be using to determine that our lift is canonical. 

\begin{theorem}[Taelman \cite{Taelman} Theorem C] \label{lenny}
Let $\O_K$ be a discrete valuation ring with perfect residue field $k$ of characteristic $p$ and fraction field $K$ of characteristic $0$. Let $X_{\O_K}$ be a projective $K3$ surface over $\O_K$ with $X_{\bar{K}}$ the geometric generic fiber and assume that $X := X_{\O_K} \otimes k$, the special fiber, is an ordinary K3 surface. Then the following are equivalent:
\begin{enumerate}
\item $X_{\O_K}$ is the base change from $W(k)$ to $\O_K$ of the canonical lift of $X$,
\item $H^2_{et}(X_{\bar{K}}, \Z_p) \cong H^0 \oplus H^1(-1) \oplus H^2(-2)$ with $H^i$ unramified $\Z_p[Gal_K]$-modules, free of rank $1, \ 20, \ 1$ over $\Z_p$ respectively. 
\end{enumerate} 
Here, the $(-1)$ and $(-2)$ denote Tate twists. 
\end{theorem}

\section{Twisted partners via Moduli space of twisted sheaves}

This section contains the main result of the article. We prove that every twisted Fourier-Mukai partner of a twisted K3 surfaces is a moduli space of stable twisted sheaves. We recall all the needed definitions and properties as well for reader's convenience and hope that it clarifies a lot of unclear details in the sparse existing literature.  

\subsection{Twisted Fourier-Mukai Partners}

Systematic study of equivalences of derived category of twisted sheaves began in thesis of Caldararu \cite{Cald}. He introduced the twisted sheaves as sheaves which glue only upto an element of the Brauer class. 
In this article, we will be studying the equivalences between the twisted derived category of K3 surfaces over an algebraically closed field of positive characteristic using another equivalent definition as sheaves on a gerbe. We refer the reader to \cite{L1} and \cite[2.1.2]{Lthesis} for an introduction to twisted sheaves as sheaves on $\G_m$-gerbes which have a particular eigenvalue for the $G_m$-action. The equivalence is shown in \cite[2.1.3]{Lthesis}.

Since the category of twisted sheaves on a variety (with a fixed Brauer class) forms an Abelian category we can form its (bounded) derived category, denoted $D^b(X, \alpha)$, where $\alpha$ is an element in $Br(X)$. 

\begin{definition}
Let $X$ be a K3 surface over a field $k$ of arbitrary characteristic with $\alpha \in Br(X)$. We define a \textbf{twisted Fourier-Mukai partner} of $X, \alpha$ as a pair consisting of a K3 surface $Y$ along with a class $\beta$ in the Brauer group of $Y$ such that $D^b(X, \alpha) \cong D^b(Y, \beta)$.    
\end{definition}

Over the field of complex numbers Yoshioka proved the following result, which describes all the twisted Fourier Mukai partners of a twisted K3 surface. We will use this result to prove a similar result over algebraically closed fields of positive characteristic. 

\begin{theorem}[\cite{Yos} Theorem 3.16, Theorem 4.3] 
Let X be a K3 surface over $\C$ with $B \in H^2(X, \Q) $ (corresponding Brauer class is denoted $\alpha_B \in Br(X)$) and $v \in \tilde{H}^{1,1}(X, B, \Z)$ a primitive vector with $<v,v>= 0$. Then there exists a moduli space $M(v)$ of stable (with respect to a generic polarizations) $\alpha_B-$twisted sheaves $E$ with $ch^B(E)\sqrt{td(X)}=v$ such that:  
\begin{enumerate}
\item Either $M(v)$ is empty or a K3 surface. The latter holds true if the degree zero part of $v$ is positive.
\item On $X' := M(v)$ one finds a $B$-field $B' \in H^2(X', \Q)$ such that there exists a universal family $\mathcal{E}$ on $X \times X'$ which is an $\alpha_B^{-1} \boxtimes \alpha_{B'}$-twisted sheaf. 
\item The twisted sheaf $\mathcal{E}$ induces a Fourier-Mukai equivalence $D^b(X, \alpha_B) \cong D^b(X', \alpha_{B'})$.
\end{enumerate}
\end{theorem}

\begin{remark}
In the case of $D^b(X) \cong D^b(X', \alpha_{B'})$, $M(v)$ is just a (coarse) moduli space of stable sheaves and in case it is fine, even $\alpha_{B'} = 0$, and we get the untwisted derived equivalence. 
\end{remark}

Moreover, there exists a Hodge theoretic  criterion for twisted derived equivalences. This allows one to use lattice theory to give the counting formula of Ma stated above. 

\begin{theorem}[\cite{HS1}, \cite{HS2}] 
Let $X$ and $X'$ be two algebraic K3 surfaces over $\C$ with rational $B$-fields $B$, respectively $B'$ inducing Brauer classes $\alpha$, respectively $\alpha'$. Then there exists a Fourier-Mukai equivalence $D^b(X, \alpha) \cong D^b(X', \alpha')$ if and only if there exists a Hodge isometry $\tilde{H}(X, B, \Z) \cong \tilde{H}(X, B', \Z)$ that respects the natural orientation of the four positive directions. 
\end{theorem}

\begin{remark}
B-fields: We will not be using them in this article. $B$-fields were introduced so that one could define Chern character for twisted sheaves in the formulation of C\u{a}ld\u{a}raru, but we replace their choice by making a choice of line bundle and working with gerbes. The Chern character then has a much more natural description which removes a lot of artificial constraints put on the definition by B -fields.
\end{remark}

\begin{remark}
Over $\C$, we have the following isomorphism \cite[Chapter 18 page 416]{HuyK3}:
$$
Br(X) \cong Hom(T(X), \Q/\Z).
$$
This implies that if we have $D^b(X) \cong D^b(Y)$ for two K3 surfaces, then this derived equivalence induces an isomorphism of the Brauer groups, $\varphi: Br(X) \cong Br(Y)$.  Moreover, the above result of Huybrechts and Stellari implies that  we have an isomorphism $D^b(X, \alpha) \cong D^b(Y, \varphi(\alpha))$. \\
In positive characteristic, we will be able to recover the first observation for ordinary K3 surfaces, since an ordinary K3 surface with a Brauer class admits a canonical lift (along with unique lifting of the Brauer class) and every Fourier-Mukai partner of it is a moduli space of stable sheaves \cite{LO}, we have that the geometric generic fibers of the canonical lifts are also derived equivalent, see \cite{Sri} for an argument in full details. Thus we have that the Brauer groups of the geometric generic fibers of the canonical lifts of $X$ and $Y$ are isomorphic but as the lifts of the Brauer classes are unique, we have that the Brauer groups of the underlying ordinary K3 surfaces are also isomorphic. 

We cannot recover the second observation by the same argument as before due to an unavailability of cohomological criterion for twisted derived equivalences.
\end{remark}

\subsection{Moduli of twisted sheaves}
Let $X$ be a K3 surface over an algebraically closed field $k$ of characteristic $p >3$, and let $\mathfrak{X} \ra X$ be a $\mathbb{G}_m$-gerbe for any integer $n$. The aim of this section is study the derived category of twisted sheaves on $\X$ by extending the results of Bragg-Lieblich\cite{BL}, Bragg  \cite{Bragg} and Yoshioka \cite{Yos} to the finite height case.

\subsubsection{Twisted Chern Character}
We begin by recalling the definition of twisted Chern character. Let $X$ be a smooth projective variety, $\alpha \in Br(X)$ a Brauer class, $\pi: \X \ra X $ be a $\mathbb{G}_m-$gerbe whose associated cohomology class is $\alpha$. Let $\alpha$ be $n$-torsion, then there exists an $n$-fold twisted invertible sheaf (use \cite[Proposition 2.1.2.5]{L1}, for the morphism $\G_m \xrightarrow{\times n} \G_m$, which induces a morphism of gerbes $\X_{\alpha} \ra \X_{n\alpha}$, then the corresponding morphism on twisted sheaves, gives us the required line bundle.), say $\mathcal{L}$, on $\X$. A choice of such a sheaf allows us to compare $n$-fold twisted sheaves on $\X$ and sheaves on $X$. 

\begin{definition} \label{cherntwisted}
Let $\mathcal{E}$ be a locally free twisted sheaf of positive rank on $\X$. The \textbf{twisted Chern character} of $\mathcal{E}$ (with respect to $\mathcal{L}$) is 
$$
ch^{\mathcal{L}}\mathcal{E} = \sqrt[n]{ch(\pi_{*}(\mathcal{E}^{\otimes n}\otimes \mathcal{L}^{\vee}))} \in A^*(X) \otimes \R,
$$ 
where by convention we choose the $n$-th root so that $rk(ch^{\mathcal{L}}(\mathcal{E})) = rk(\mathcal{E})$. Moreover, the \textbf{Mukai vector} of $\E$  is
$$
v^{\L}(\E) = ch^{\L}(\E)\sqrt{Td(X)},
$$
where $Td(X)$ is the Todd genus of $X$.   
\end{definition}

In the tame setting, our $\G_m$-gerbe, can be "replaced" by a $\mu_n$-gerbe, which being tame is  a Deligne-Mumford stack with $X$ as its coarse moduli space. Moreover, since for a K3 surface, cohomological Brauer group is the same as (Azumaya) Brauer group, we have from \cite[Theorem 3.6]{EHKV}, that any $\mu_n$-gerbe is a quotient stack. Then the conditions for the result of Kresch \cite[Proposition 5.1]{KreschDM} are satisfied and our gerbe admits resolution property. Hence any twisted sheaf on $\X$ admits a resolution by locally free twisted sheaves of positive rank, thereby allowing us to extend the definition of twisted Chern character and twisted Mukai vector, by additivity, to all of Grothendieck group of twisted sheaves, $K(X, \alpha)$. The following very easy lemma gives us the image range of the twisted Chern character. 

\begin{lemma}
Twisted Chern character actually has image in $A^*(X) \otimes \Q$. Thus, also in $A^*(X) \otimes \Q_l$. 
\end{lemma}

\subsubsection{Lattices associated to a K3 surface}
\begin{definition}
Let $X$ be a K3 surface over an algebraically closed field $k$. Let $N(X)$ be the N\'eron-Severi lattice of $X$, we define the \emph{extended N\'eron Severi lattice of $X$} to be free $\Z$-module
$$
\tilde{N}(X) = \langle (1,0,0) \oplus N(X) \oplus (0,0,1) \rangle
$$
equipped with the Mukai pairing
\begin{eqnarray*}
\tilde{N}(X) \otimes \tilde{N}(X) \longrightarrow &\Z \\ 
\langle (a, b, c) , (a', b', c') \rangle \mapsto &-a.c' + b.b' + a'.c . 
\end{eqnarray*}
\end{definition}

\subsubsection{Height is a Twisted Derived Invariant} \label{HeightTDI}
Let $\X \ra X$ and $\mathcal{Y} \ra Y$ be $\G_m$-gerbes over K3 surfaces $X$ and $Y$ respectively. The following is the Grothendieck-Riemann Roch formula for the projection $\pi_{\X}: \X \times \mathcal{Y}   \ra \X$, with $\pi_X: X \times Y \ra X$ the induced map on coarse moduli spaces.  
\begin{lemma} \label{TwistedGRR} 
For any $\alpha \in K^{(1,0)}(\X \times \mathcal{Y})$, we have 
$$
ch_{\X} (\pi_{\X*} \alpha) = \pi_{X*}(ch_{\X \times \mathcal{Y}}(\alpha).Td(\pi_X)).$$
\end{lemma}

\begin{proof}
We will extend the proof of \cite[Lemma 4.1.4]{BL}. The only place where we need a non-obvious generalization is to find a trivialization  of $\X$ by a finite flat cover (this was provided by absolute frobenius in case of $\mu_p$-gerbes). For this, use \cite[Proposition 3.3.2.6]{L1} and \cite[Lemma 2.3.4.2]{L1} along with \cite[12.3.10]{OlssonAS} as for K3 surfaces, the cohomological Brauer group and scheme theoretic Brauer group are isomorphic. 
\end{proof}

Height of a K3 surface was defined above in definition \ref{heightdef}. Now we show that height of a K3 surface is a twisted derived invariant. More precisely, 

\begin{lemma}(cf. \cite[Proposition 4.1.7]{BL})
 Let $\X \ra X$ and $\mathcal{Y} \ra Y$ be two $\G_m$-gerbes over K3 surfaces $X$ and $Y$ respectively, and $\mathcal{P} \in D^{(1,1)}(\X \times \mathcal{Y})$ a perfect complex of twisted sheaves inducing an equivalence of derived categories $D^{(-1)}(\X) \cong D^{(1)}(\mathcal{Y})$ (i.e., $\X$ and $\mathcal{Y}$ are derived equivalent), then the height of $X$ = height of $Y$.  
\end{lemma}

\begin{proof}
The follows from descent of a twisted derived equivalence to an isometry of the F-isocrystals. We will follow the proof in \cite[Proposition 4.1.7]{BL}. 
To show that 
$$
\Phi_{v_{\X \times \mathcal{Y}}(\mathcal{P})}^{crys}: \tilde{H}(X/K) \ra \tilde{H}(Y/K)
$$
preserves the F-isocrystal structure, we have to show that the following three structures are preserved:

\begin{enumerate}
    \item Isomorphism of K-vector spaces: This follows from the easy observation that $\Phi_{\O_{\Delta}} = id $ on the derived categories and hence so on the cohomology.  
    \item Isometry: Same as in untwisted case. 
    \item Compatibility with Frobenius: For this one has to just show that for $ch(\F)  \in H^*_{crys}(X/K)$, $$
    q_*(p^*(F_X^*(ch(\F)).v(\mathcal{P}))= F_Y^*q_*(p*(ch(\F)).v(P))
    $$
    which follows easily from the  K\"unneth formula compatibility with absolute Frobenius endomorphism. Here $p$ (resp. $q$) is the projection from $X \times Y$ to $X$ (resp. $Y$). This is same as in untwisted case. 
\end{enumerate}
\end{proof}

\begin{lemma}\cite[Lemma 3.3.7]{LMS}
Twisted Chern character for a K3 surface is $l$- adic integral. 
\end{lemma}

\subsubsection{Stability of twisted sheaves}

Fix a polarization $H$ on $X$. Recall the following definitions
\begin{definition}
If $\E$ is a $\X$-twisted sheaf, then the \textbf{geometric twisted Hilbert polynomial} of $\E$ is the function 
$$
P_{\E}(m) =deg\big(ch^{\mathcal{L}}_{\X}(\E(m)).{Td(X)}\big),
$$
where $\E(m) = \E \otimes p^*(H)^{\otimes m}$. Moreover, the \textbf{reduced twisted Hilbert polynomial} is defined as 
$$
p_{\E}(m)=\frac{1}{a_d} P_{\E}(m), 
$$
where $a_d$ is the leading coefficient of $P_{\E}$. 
\end{definition}
\begin{remark}
Note that the choice of an $n$-fold twisted invertible sheaf $\L$ on $\X$ is essentially the same as a choice of preimage of $\alpha$ under the map $H^2(x, \mu_n) \twoheadrightarrow Br(X)[n]$, which is giving a $\mu_n$-gerbe $\X'$ along with an isomorphism of its associated $\mathbb{G}_m$-gerbe with $\X$. 
\end{remark}
\begin{definition}
An $\X$-twisted sheaf $\E$ is \textbf{stable} (resp. \textbf{semistable}) if it is pure and for all proper non-trivial subsheaves $\F \subset \E$ 
$$
p_{\F}(m) < p_{\E}(m) \ ( \text{resp.} \ \leq) \ \forall m >> 0.
$$
\end{definition}

\begin{definition}
Let $X \ra S$ be a realtive K3 surface, $\X \ra X$ a $\mathbb{G}_m$-gerbe, $v$ a global section of $\tilde{H}(X/S) \otimes \Q$ and $H$ a relative polarization of $X/S$. The \textbf{moduli space of $\X$-twisted stable sheaves with twisted Mukai vector $v$} is the stack $\mathcal{M}_{\X/S}(v)$ on $S$ whose objects over an $S$-scheme $T$ are $T$-flat $\X_T-$twisted sheaves $\E$ locally of finite presentation such that for each geometric point $t \in T$, the fiber $\E_t$ is $H_t$-stable and has twisted Mukai vector $v_t$.   
\end{definition}

The stack $\mathcal{M}_{\X/S}$ is an algebraic stack for the case that the $\mathbb{G}_m$-gerbe $\X/S$ comes from a $\mu_n$-gerbe, where $n$ is coprime to characteristic of $X$ as shown in \cite[Proposition 2.3.1.1, 2.3.2.11]{L1}. Moreover, the stack $\mathcal{M}_{\X/S}(v)$ is
an algebraic substack of finite type over $S$.

For reader's convenience, we have made a comparison table, it is a good way to see where things differ and where not:

Let $X$ be a smooth projective variety $k =\bar{k}$, $\X \ra X$ a $\mu_n$-gerbe $n$ invertible or $\G_m$-gerbe.
\begin{table}[hbt!]
\begin{tabular}{|l | c | c| }
\hline
Property & Coherent sheaf & Coherent Twisted sheaf \\
\hline \hline
Support &$ \{x \in X | \F_x \neq 0 \}$  & $\{x  \in |\X| | x^*\F \neq 0 \}$ \\
	& closed set & closed\\
 &&\\
Euler char & HRR* & $\chi(\F) := [\mathcal{I}(\X):\X]deg(ch(\F) Td_{\X})$ \\
&& \\
Hilbert poly & $n \mapsto P_{\F}(n) = \chi(\F(n))$ & $n \mapsto P_{\F}(n) = \chi(\F(n))$ \\
&Integral coeff. $\alpha_i$ & Rational coeff. $\alpha_i$ \\
&&\\
Rank ($d = \dim \X$)  & $rk(\F) = \frac{\alpha_d(\F)}{\alpha_d(\O_X)}$ & $rk(\F) = \frac{\alpha_d(\F)}{\alpha_d(\O_{\X})}$ \\
Else $0$ &$d = \dim Supp(\F)$ & $d = \dim Supp(\F)$ \\
&& \\
Degree $deg(\F)$ & $\alpha_{d-1}(\F)$ &$\alpha_{d-1}(\F)$ \\
$d = \dim \X$ &$- rk(\F).\alpha_{d-1}(\O_{X})$ & $- rk(\F).\alpha_{d-1}(\O_{\X})$\\
&& \\
Semistable & For any $\mcG \subset \F$  & For any $\mcG \subset \F$  \\
& $\alpha_d(\F)(P_{\mcG})\leq  \alpha_d(\mcG)P_{\F}$& $\alpha_d(\F)(P_{\mcG}) \leq  \alpha_d(\mcG)P_{\F}$\\
&& \\
Stable & For any $\mcG \subset \F$ & For any $\mcG \subset \F$  \\
& $\alpha_d(\F)(P_{\mcG}) <  \alpha_d(\mcG)P_{\F}$& $\alpha_d(\F)(P_{\G}) <  \alpha_d(\mcG)P_{\F}$\\
&& \\
Slope & $\mu(\F) = \frac{deg \F}{rank \F}$ & $\mu(\F) = \frac{deg \F}{rank \F}$ \\
$d= dim(\X)$ & & \\
&&\\
$\mu$-semistable & $\F$  pure, for any $\mcG \subset \F$ & $\F$  pure, for any $\mcG \subset \F$  \\
$d= dim \X$ & $\mu(\mcG) \leq  \mu(\F)$& $\mu(\mcG) \leq  \mu(\F)$\\
&&\\
$\mu$-stable & $\F$  pure, for any $\mcG \subset \F$ & $\F$  pure, for any $\mcG \subset \F$  \\
$d= dim \X$ & $\mu(\mcG) <  \mu(\F)$& $\mu(\mcG) < \mu(\F)$\\
\hline
\end{tabular}
\end{table}

*HRR: Hirzebruch Riemann-Roch. In case of coherent sheaves we define Euler characteristic in terms of alternate sum of dimensions of coherent cohomology groups and Hirzebruch Riemann Roch then implies that we can compute it using the chern characters. On the other hand for twisted sheaves we use the chern characters to define the Euler characteristic. 

\begin{remark}
Don't divide by zero: rank, deg and $\mu$-(semi)stablity are only defined for sheaves with dimension of support equal to the dimension of $X$ (resp. dim $\X$), otherwise one would be dividing by zero as $\alpha_d(\O_X)$ can be zero for non-maximal d. This happens, for example, for d=1 and X = K3 surface.  
\end{remark}

\begin{definition}
A polarization $H$ is said to be \emph{$v$-generic} if any $H$-semistable twisted sheaf with Mukai vector $v$ is $H$-stable. 
\end{definition}

\begin{definition}
Primitive Mukai vector in (extended) Neron-Severi lattice is an element which is not an integral multiple of another element in  the lattice. 
\end{definition}

\textbf{Numerical Criterion:}
Note that if we have a locally free $\mcX$ twisted sheaf, $\E$ on $\X \ra X$ with twisted Mukai vector $v$, then its Mukai paring with $v(k(x))$ for a skyscraper sheaf is:
$$
<v, (0,0,1)> = rk(\E) 
$$
and since $rk(E)$ is divisible by $\text{ind}(\X)= n$, where $n$ is the order of class of $\X$ in $Br(X)$, we have that $<v, v(k(x))>$ is divisible by $n$.  

\begin{lemma} (\cite[Proposition 4.1.18]{BL})
Let $\mathcal{X}\ra X$ be a $\mu_n$ gerbe on a proper smooth family of  K3 surfaces over a Henselian DVR $R$, and let $H$ and $v$ be as in Definition (of moduli spaces). Suppose that $v$ restricts to a primitive  element of the twisted Neron-Severi group of each geometric fiber. If $H$ is $v$-generic in each geometric fiber of $\mathcal{X}/R$, then the moduli space $M_{\mathcal{X}}(v)$ of $H$-stable twisted sheaves on $\mathcal{X}$ with twisted Mukai vector $v$ is either empty or a $\G_m$-gerbe over a proper smooth scheme over $R$. 
In particular, every fiber of $M_{\mathcal{X}}(v) \ra Spec (R)$ is a $\G_m$-gerbe over a K3 surface if and only if one geometric fiber is a $\G_m$-gerbe over a K3 surface. 
\end{lemma}

\subsubsection{Lifting a $v$-generic polarization} \label{liftingv-generic}

Let $X$ be a K3 surface over an algebraically closed field $k$ of characteristic $p > 0$. Assume we have a lift of $X$, denoted $\mcX_A$, to a DVR $A$ with residue field $k$ and fraction field $K$. 
Let $v \in A^*(X)$ be a Mukai vector, assume that the above lift is such that $v$ lies in the image of the specialization map on the Chow groups \cite[Corollary 20.3]{FultonIT}. This condition in the untwisted case is ensured by taking our Mukai vector to be of the form $v= (r, l, s)$, where $l$ is the Chern class of an ample line bundle. Then choosing a lift to which this line bundle lifts.   

This trick does not look like it will work in twisted case-The issue being that for untwisted sheaves our Mukai vector of the sheaf will return just zero! 

We will abuse the notation and pick a Mukai vector in the preimage of $v$ under the specialization map of the chow group and denote it still by $v$ (both for generic fiber and the lift over the DVR). We would like to show that if a polarization, say $H$, on $X$ was $v$-generic, then the lift of it, say $H_W$ (resp. $H_K$) to $\mcX_W$ (resp. $\mcX_K$) will be also be $v$-generic.

To ensure that we will be able to do it, we would need a special type of lifts of a K3 surface, only possible for K3 surfaces of finite height, the Picard preserving lifts as constructed by Lieblich-Maulik in \cite{LM}. Let us recall walls of type $(r, \Delta)$ from \cite[Section 4.C]{HL} and \cite[Section 10.2.4]{HuyK3} (c.f. \cite[Section 3.2.1]{L1}). They are defined case wise depending on the dimension of the support of sheaves. We can do this since fixing the  Mukai vector of sheaves fixes their Hilbert polynomial and thereby fixing the dimension of their support. 

\textbf{Case 0: $\dim (\Supp \F) =0$.}
Note that, if the dimension of support of a coherent sheaf is zero, then it is supported at only finitely many points, and it is reduced Hilbert polynomial is always 1. Thus, it is stable if and only if it contains no proper subsheaves, i.e., it is supported only on one point and it has length 1 as a sheaf. Moreover, note that a change of polarization has no effect on stability of such a sheaf, since tensoring a skyscraper sheaf with a line bundle gives us back the skyscraper itself.

\textbf{Case 1: $\dim (\Supp \F) =1$}. We will follow the definition in \cite[Section 1.4]{YoshiokaAS}. Recall that by the definition of the first Chern class of a semistable (hence, pure) dimension one sheaf that it is the class of the  dimension 1  closed subscheme (possibly non-integral) of $X$ defined by the zeroth Fitting ideal of $\F$, see \cite[Tag 0C3C]{SP} for the definition. This is an effective divisor of $X$. Let  $H$ be an effective divisor such that $c_1(H) = c_1(\F)$ and $(H)^2 >0$.   
Semi-stability of $F$ with respect to any polarization $L$ can be phrased as, 
$$
\frac{\chi(\E)}{(c_1(\E), L)} \leq \frac{\chi(\F)}{(c_1(\F), L)}, 
$$
where $E \neq 0$ is a proper subsheaf of $F$. Stability is equivalent to strict inequality. 
Assuming $\chi(F) \neq 0$, we define the wall and chamber structure on $Num(X)$ for dimension 1 sheaves. 

\begin{remark}
This assumption is easy to satisfy, just tensor by the line bundle corresponding to $H$.
\end{remark}
 
Let $E$ be a subsheaf of $\F$, define $\xi := \chi(E)c_1(\F) -\chi(\F)c_1(E)$. Now $\chi(\F) \neq 0$ implies that if $c_1(E) \not\in \Q c_1(\F)$, then $\xi \neq 0$ and for such a  $\xi \neq 0$, we define the \textbf{walls} $W_{\xi}$ as follows:
$$
W_{\xi} := \{ x \in Amp(X) | (x, \xi) = 0\} \subset Amp(X). 
$$
We call the connected components of $Amp(X) \setminus \bigcup_{\xi} W_{\xi}$ a \textbf{chamber}.

Then from  \cite[Lemma 1.2]{YoshiokaAS} and the discussion above it, we have that the number of walls are finite and for a generic polarization we have that all semistable sheaves are stable.   

\textbf{Twisted case:} Note that any twisted dim 1 sheaf on $\X$, a gerbe over $X$ a K3 surface is actually untwisted. Indeed, Any dimension 1 sheaf on $\X$ can be seen as a torsion free on a dimension 1 closed substack of $\X$, denoted $\X_C \hookrightarrow \X$. Note that the coarse moduli space of $\X_C$ is a  closed dimension 1 subscheme of $X$. (First note that properness of all the involved stacks, imply that the coarse moduli space of $\X_C$ is a dimension 1 scheme and then pushforward the structure sheaves to get that the coarse moduli space is indeed defined by the pushforward of the ideal sheaf using that pushforward is an isomorphism for untwisted sheaves.) Moreover, $\X_C$ is a gerbe over its coarse moduli space, denoted $C$, as we have a Cartesian square
\begin{equation*}
   \xymatrix{
    \X_C  \ar@{^{(}->}[r] \ar[d] & \X \ar[d] \\
C \ar@{^{(}->}[r]& X.
}
\end{equation*}

Then from Giraud Section V.1 or \cite[Theorem 3.11]{BSdecomp}, we see that the class of $\X_C$ over $C$ in $H^2(C, \G_m)$, is the pullback of the class of $\X \ra X$ in $H^2(X, \G_m)$ under the imbedding $C \hookrightarrow X$. And thus from \cite[Proposition 4.9]{BSdecomp}, we have that the pullback preserves the $\chi$-twisted sheaves decomposition, and hence we have a 1-twisted sheaf on $\X_C$. 

But recall that for a dimension 1 scheme over an algebraically closed field, the Brauer group is trivial. Indeed, from \cite[Proposition 7.2.1, Remark 7.2.2]{CTS}, we have an injective morphism 
$$
Br(C) \hookrightarrow Br(\tilde{C}) \oplus_{x \ closed} Br(k(x)),
$$ 
where $\tilde{C}$ is the normalization of $C$ and by Tsen's theorem, its Brauer group is zero and as we are over algebraically closed fields $Br(k(x)) =0$. Thus, $Br(C) =0$. 

And hence we have $\X_C$ is the trivial gerbe and the 1-twisted sheaf we started with actually comes from an untwisted sheaf on $C$.

\textbf{Case 2: $\dim (\Supp \F) =2$}. For this we use the definition in \cite[Section 4.C]{HL} and for the twisted case see \cite[Proposition 4.1.14]{BL}.

From the above discussion of walls and $v$-generic polarizations, it is clear that if we find a lift our K3 surface such that  the natural specialization map of Picard lattice: 
$$
sp: Pic(\mcX_{\bar{K}}) \ra Pic(X)
$$
is an isomorphism and it preserves the ample cone, then every $v$-generic polarization lifts to a $v$-generic polarization.
A Picard preserving lift can be constructed for any K3 surface of finite height as in \cite[Corollary 4.2]{LM}. Also recall that canonical lift of an ordinary K3 surface is a Picard preserving lift.   Moreover, by \cite[Corollary 2.4]{LM}
for any Picard preserving lift, the specialization map also preserves the ample cone.

Thus, we can now use lifting to char 0 arguments to prove the theorems below.

\begin{theorem}(cf. \cite[Theorem 4.1.19]{BL}) \label{moduli}
Let $\X \ra X$ be a $\mathbb{G}_m$-gerbe on a K3 surface of height $h$, where the order of class of the gerbe in $Br(X)$ is coprime to $p$  and $v = (r, l, s)$ a primitive Mukai vector with $v^2=0$. If $H$ is $v$-generic, then the moduli stack $\mathcal{M}_{\X}(v)$ of $H$-stable twisted sheaves on $\X$ with twisted Mukai vector $v$ is either empty or satisfies
\begin{enumerate}
\item $\mathcal{M}_{\X}(v)$ is a $\mathbb{G}_m$-gerbe over a K3 surface $M_{\X}(v)$ of height $h$.
\item the universal sheaf $\mathcal{P}$ on $\mathcal{M}_{\X}(v) \times \X$ induces a Fourier-Mukai equivalence
$$
\Phi_{\mathcal{P}}: D^{(-1)}(\mathcal{M}_{\X}(v)) \ra D^{(1)}(\X),
$$ 
and
\item the $\mathbb{G}_m$-gerbe is trivial if and only of there exists a vector $w \in A^*(X) \otimes \Q$ such that $v.w$ is co-prime to n. 
\end{enumerate}
Moreover, the moduli space $\mathcal{M}_{\X}(v)$ is non-empty if one of the following is true:
\begin{enumerate}[label=(\roman*)]
    \item the degree zero part of $v$ is positive, i.e., $r > 0$. 
    \item In case $r = 0$, we have $l$ is class of an effective divisor and $s \neq 0$. 
    \item In case $v =(0,0,1)$, in which the moduli space is isomorphic to $X$.
    \end{enumerate}

\end{theorem}

\begin{proof}

\textbf{Case $r > 0$:} We remark that the above result is just \cite{LMS}[Proposition 3.4.2], in the special case of $\text{rk}(v) = n$. Note that the from the period index theorem for surfaces, the $\text{rk}(v)$ is a multiple of $n$. And the above case is enough for our purposes due to the following:

\textbf{Claim:} Given a rank $r$ twisted sheaf on a $\mu_n-$gerbe $\X \ra X$, we can construct a $\mu_r-$ gerbe $\X_r$ over $X$, with the same Brauer class as $\X$. Moreover, the natural map $\X_r \ra \X$ serves to identify the stack of semistable sheaves via pullback.

 This follows from construction \cite[2.2.2.11]{Lthesis} using the Azumaya algebra $\mathcal{E}nd(\F)$ of degree $r$, where $\F$ is the rank $r$ locally free twisted sheaf. The class of the constructed gerbe is the same as the class of the Azumaya algebra which is the same as the class of the gerbe on which $\F$ is a 1-twisted sheaf.

Note that the non-emptiness of $\mathcal{M}$ will follow just from lifting our K3 surface to characteristic $0$.

\textbf{Case $r =0$}: These include both the cases (ii) and (iii) above.  We will prove  that the moduli space is actually isomorphic to moduli of untwisted sheaves in characteristic zero and then lift the char p moduli space  to char 0 and conclude the same for them. 

From the discussion the above the theorem we see that the closed points of the moduli space of rank 0 twisted sheaves actually corresponds to untwisted sheaves.

Moreover, using the boundedness of a set of twisted sheaves with fixed Mukai vector, we can actually untwist a family sheaf of twisted sheaves to get a family of untwisted sheaves. Indeed, for example, in case of dimension 1 sheaves (dimension 0 sheaves is similar and even easier), let $\F$ be a family of twisted sheaves defined on $\X \times S$, i.e. $\F \in \mathcal{M}_{\X}(v)(S)$, for $S$ a $k$-scheme. Now, for each point $s \in S$ the sheaf $\F_s$ can be seen as a twisted sheaf on $\X$, supported on a  (gerby) curve $\X_{C_s}$ and as discussed above $\X_{C_s}$ is actually a trivial gerbe over a curve and hence has a twisted line bundle, whose dual will untwist the sheaf $\F_s$. We denote the line bundle by $\L_s$ and its dual by $\L^{\vee}_s$. Now, the pushforward of $\L^{\vee}_s$ to $\X$ for each $s \in S$ gives a set of dim 1 twisted  sheaves on $\X$, but such a set is bounded for example  by \cite[Theorem 4.12]{NironiSS}. And we can get a family of sheaves on $\X \times S$, such that its fibers are $\L^{\vee}_s$. Now we tensor this family, denoted $\L$ with $\F$, to get a family of untwisted dim 1 sheaves on $\X$ (and hence sheaves on $X$) and thus we can define a functor (actually a natural transformation) between the moduli stacks as follows:
\begin{eqnarray*}
   \varphi: \mathcal{M}_{\X}(v)(S) \ra &\mathcal{M}_X(v)(S) \\
    \F \mapsto & \F \otimes \L.
\end{eqnarray*}
(The definition of the morphism for $Hom$ is also similar.)

This map induces a map at the level of coarse moduli spaces. Indeed, the composition with the coarse moduli space morphism for the moduli of untwisted sheaves $\mathcal{M}_{\X}(v)\xrightarrow{\varphi}\mathcal{M}_{X}(v) \ra M_X(v)$ factors via the coarse moduli space of twisted sheaves because of the universal property of coarse moduli spaces of stacks. 

Now, note that at the morphism of coarse moduli space is bijective at the level of closed points and hence when the characteristic of the algebraically closed ground field is zero, we actually have an isomorphism of the coarse moduli spaces and hence we have an isomorphism of the moduli stacks as well as they are gerbes and any morphism between $\mu_n$-gerbes (or even as $\G_m$-gerbes) over isomorphic bases is an isomorphism. 

In case the ground field is an algebraically closed field of characteristic $p$, we lift the K3 surfaces of finite height to characteristic zero using a Picard preserving lift (see the discussion above the theorem) and hence we have we can lift our $v$-generic polarization as a $v$-generic polarization as well. Now consider the lifts and the moduli spaces of (twisted) sheaves over them. Let us fix the notation: $\mcX_W$ is the lift of $X$ with (geometric) generic fiber $\mcX_{\bar{K}}$, the corresponding gerbes over them $\X_W$ and $\X_{\bar{K}}$, the moduli spaces $\mathcal{M}_{\X_W}(v)$ (resp. $\mathcal{M}_{\X_W}(v)$ ) over $\X_W$ (resp. $\X_{\bar{K}}$) of twisted sheaves. These stacks(gerbes) are smooth and proper  and hence also their coarse moduli spaces.  Now from the characteristic zero result that we just proved before, we know that at the level of generic fibers the stack of twisted sheaf is just the stack of untwisted sheaves. Also the coarse moduli spaces are isomorphic. Then using the fact that a smooth proper surface over an algebraically closed field is projective, we can apply Matsusaka and Mumford \cite[Corollary 1]{MM}, we get that the coarse moduli spaces of special fibers are K3 surfaces, when not empty and also that the isomorphism can be extended to an isomorphism of the relative coarse moduli space. And again as the stacks we are considering are gerbes we get an isomorphism of the stacks as well. \end{proof}

\subsection{Twisted Fourier-Mukai partners}

We now prove that every twisted Fourier-Mukai partner of a twisted K3 surface is isomorphic to a moduli space of twisted sheaves on it. More precisely, this can be  stated as follows.

\begin{theorem} \label{TwistedpartnersModuli}
Let $\X \ra X$ (resp. $\mathcal{Y} \ra Y$) be $\G_m$-gerbe over $X$ (resp. over $Y$) and $\mathcal{P} \in D^{(1,1)}(\X \times \mathcal{Y})$ be a twisted complex inducing a derived equivalence 
$$
\Phi_{\mathcal{P}}: D^{(-1)}(\mathcal{Y}) \xrightarrow{\cong} D^{(1)}(\mathcal{X}),
$$
then $Y$ is isomorphic to $M_{\X}(v)$ and $\mathcal{Y}$ is a isomorphic as a $\G_m$-gerbe to $\mathcal{M}_{\X}(v)$, where $M_{\X}(v)$ is the moduli space of stable twisted sheaves with Mukai vector $v = v(\mathcal{P}|_{\X \times y})$ on $\X$, for some fiber (hence all fibers).  
\end{theorem}

We will be prove this theorem in the same way its untwisted cousin was proven in \cite{LO}. To recall, we will first show that every filtered twisted derived equivalence which preserves the ample cone admits a lifting to characteristic zero and there it induces an isomorphism on the corresponding geometric generic fibers of the lifted K3 surfaces and hence by Matsusaka-Mumford induces an isomorphism on the special fibers. That is, we would have show that filtered twisted derived equivalences induce an isomorphism of the underlying twisted K3 surfaces. Then for any twisted Fourier-Mukai partner of a twisted K3 surface we compose with the inverse of the twisted derived equivalence induced by the universal sheaf of the moduli space of twisted sheaves to get a filtered twisted equivalence  and hence we can conclude that the twisted Fourier-Mukai partner is just a moduli space of stable twisted sheaves. 

Let us start with the definition of filtered derived equivalences. 

\begin{definition}
Let $\X$ (resp. $\Y$) be a $\G_m$-gerbe over $X$ (resp. $Y$) a K3 surface, then we say a twisted derived equivalence 
$$
\Phi_{P}: D^{(-1)}(\X) \ra D^{(1)}(\Y)
$$
is \textbf{filtered} if and only if the induced isometry at the level of Chow groups sends $(0,0,1)$ to $(0,0,1)$, i.e.
\begin{eqnarray*}
\Phi_{P}^A: A^*(X) &\ra A^*(Y) \\
(0,0,1) &\mapsto (0,0,1). 
\end{eqnarray*}
\end{definition}

\begin{lemma}[cf. \cite{BraggYang}]
If $\Phi_P: D^{(-1)}(\X) \ra D^{(1)}(\Y)$ is a filtered derived equivalence of $\mu_n$-gerbe over $X$ with a  $\mu_m$-gerbe over $Y$, then $m =n$. That is, they have same orders. 
\end{lemma}

\begin{proof}
We will show that $m | n$ and the converse relation  can be checked using the inverse equivalence of $\Phi_P$. 
Recall that we have a locally free 1-twisted sheaf of rank $n$ on $\X$, denoted $\E$ and taking its determinant gives us an n-twisted line bundle on $\X$. Use the dual of this line bundle (resp. its n-tensor power) to define the Chern characters of 1-twisted sheaves (resp. n-twisted sheaves) as in definition \ref{cherntwisted}.Hence, $ch^{det(\E)}(\E) = (n, 0,s)$. Now consider the image of this sheaf under $\Phi_P$, i.e, $\Phi_P(\E)$. This is a perfect complex of some rank $n'$ and taking its determinant gives us a line bundle, which is $n'$-twisted. For $\Y$ we will take Chern characters with respect to the line bundle $det(\Phi_P(\E))$ and hence we get $ch^{det(\Phi_P(\E))}(\Phi_P(\E)) = (n', 0, s')$. Now since $(0,0,1)$ is mapped to $(0,0,1)$ by $\Phi_P$, and $\Phi_P$ is an isometry we have 
\begin{eqnarray*}
n &= <(n,0,s), (0,0,1)> = <v(\E), (0,0,1) >  \\ &=<v(\Phi_P(\E)),(0,0,1)>  = <(n',0,s), (0,0,1)> = n'.
\end{eqnarray*}

Thus we get a rank n-twisted perfect complex on $\Y$, whose determinant will be a $n$-twisted invertible sheaf on $\Y$. It can be easily seen that existence of such a sheaf implies that the order of the gerbe $\Y$ is divisible by $n$. On doing the same argument for the inverse derived equivalence gives us that $n=m$.  
\end{proof}

\begin{lemma}
We can assume that a filtered derived equivalence sends $(1,0,0)$ to $(1,0,0)$ in the Chow groups. 
\end{lemma}

\begin{proof}
From the previous lemma, we know that the gerbes $\X$ and $\Y$ have same orders and we will be making the same choice of line bundles for computing the chern characters. The claimed result will follow from the following easy computation: 

Using $\Phi_P^A$ is an isometry, we get that $<(n,0,s').(n,0,s')> =-2ns$, thus we see that $s =s'$. Thus $(n,0,s)$ is mapped to $(n,0,s)$ and now since $(0,0,1)$ is mapped to $(0,0,1)$ and the morphism at the level of chow groups is a group homomorphism we get that $(1,0,0)$ is mapped to $(1,0,0)$.
\end{proof}


\begin{remark}[Spherical Twists] 
Unlike the case of untwisted K3 surfaces, inside the twisted derived category of a K3 surface, it can happen that there are no spherical object. Recall that a spherical object is an object $E$ of the (twisted) derived category $D^b(X, \alpha) = D^b_{(1)}(\X)$ such that 
\begin{eqnarray*}
Ext^i(E, E)= 
\begin{cases}
      1 & i = 0, 2 \\
     0 & \text{otherwise}
    \end{cases}   
\end{eqnarray*}
In the untwisted case (i.e, $\alpha =0$), we always have $\O_X$ or even $\L$, a line bundle, all give spherical objects. But these objects cannot be twisted to give a twisted spherical object unless the gerbe is trivial. Even more is true, namely that there exist twisted K3 surfaces which do any admit any spherical object, for example, see \cite[Lemma 3.22]{HMSgenericK3}.  Thus, we don't have the corresponding spherical twists for them giving a twisted derived autoequivalence, although we do have twisted derived autoequivalences induced by tensoring every twisted complex by a line bundle.  

However, in case we have (-2)-curves, $C$ (e.g. any smooth rational curve) in $X$, we have sheaves $\O_C$, and as discussed above in Case 1 of Section \ref{liftingv-generic}, we can consider them as twisted sheaves on $\X$ and these give spherical objects (twisted as well as well untwisted) and the corresponding spherical twist \cite[Definition 8.3]{HuyFM} is gives a (twisted) derived autoequivalence. The proof is same in twisted case as in untwisted, see \cite[Proposition 8.3]{HuyFM}. As an example of K3 surface with spherical twists, let us mention that existence of non-ample line bundles on a K3 surface implies that there exist (-2) curves on $X$ \cite[Corollary 8.1.6]{HuyK3} and hence such a K3 surface admits (twisted) spherical twists.
\end{remark}

\begin{proposition}[c.f \cite{LO} Theorem 6.3] \label{LOlifting}
Let $(X, \alpha)$ and $(Y, \beta)$  be two twisted K3 surfaces over an algebraically closed field $k$, with $\X$ and $\Y$ the corresponding gerbes respectively. Let $P \in D^b_{(1,1)}( \Y \times \X)$ be a perfect complex inducing a twisted derived equivalence $\Phi_P: D^b_{(-1)}(\Y) \ra D^b_{(1)}(\X)$ on the derived categories. Assume that the induced map on the Chow group satisfies:
\begin{enumerate}
\item $\Phi(1, 0, 0) = (1, 0, 0)$, 
\item the induced isometry $\kappa: Pic(Y) \rightarrow Pic(X)$ sends $C_Y$, the ample cone of Y, isomorphically to either $C_X$ or $-C_X$, the $(-)$ample cone of $X$.   
\end{enumerate}  
Then there exists an isomorphism of infinitesimal deformation functors $\delta: Def_{\X} \rightarrow Def_{\Y}$ such that 
\begin{enumerate}
\item $\delta^{-1}(Def_{(\Y,L)}) = Def_{(\X, \Phi(L))}$;
\item for each augmented Artinian $W$-algebra $W \rightarrow A$ and each $(\X_A \rightarrow A) \in Def_{\X}(A)$, there is an object $P_A \in D^b(  \delta(\X_A) \times_A \X_A)$ reducing to $P$ on $\Y \times \X$. 
\end{enumerate}
\end{proposition}

\begin{remark}
Before we begin the proof, we give a remark that the conditions of the theorem are satisfied by any filtered twisted derived equivalence: (1) from the lemma above and for (2), the same proof as in \cite[Lemma 6.2]{LO} works using the remark above.
\end{remark}

\textbf{Idea of the proof:} Given a deformation of $\X$ we want to find a deformation of $\Y$ such that the Fourier-Mukai kernel $P$ can also be lifted along with. This problem can be reformulated using moduli spaces of objects in the derived categories as follows, consider the moduli stack(s) of (twisted) universally gluable relatively perfect complexes of  $\X$-twisted sheaves (resp. $\X_A$-twisted), denoted $Tw_{\X/k}$ (resp. $Tw_{\X_A/A}$) then $P$ gives a object $Tw_{\X /k}(\Y) \hookrightarrow Tw_{\X_A/A}(\Y)$, corresponding to a morphism $\Y \ra Tw_{\X_A/A} $, now existence of a morphism $\Y_A \ra Tw_{\X_A/A}$ we do the job for us. The existence of such a morphism will follow from the smoothness of the moduli stacks and their coarse moduli spaces. This is just a generalization of the proof of \cite[Theorem 6.3]{LO} for the twisted setting following the twisted analogue's developed in \cite[Section 3, Section 4]{Re}. The proof is added for reader's convenience. 

Let us also remark that the condition of universally gluable is the condition which allows us to glue morphisms in the derived categories and get our moduli stack, without this condition one doesn't have a stack but if one is willing to work with the more general $\infty$-category of perfect complex (whose homotopy category is the derived category of coherent sheaves), we can remove the gluable condition and construct a derived moduli stack, see To\"en-Vaquie \cite{TV}.

\begin{proof}
Let $sTw_{\X/k}(0,0,1)^0$ (resp. $sTw_{\X_A/A}(0,0,1)^0$) be the coarse moduli space of the moduli stack of universally gluable simple twisted complexes on $\X$ (resp. $\X_A$) with trivial determinant (see \cite[Definition 3.1, 4.1]{Re}) with Mukai vector $(0,0,1)$. 

\textbf{Claim:} The moduli space $sTw_{\X_A/A}(0,0,1)^0$ is smooth. 

For this, like in proof of \cite[Proposition 4.2]{Re}, we use \cite[Tag 0APP]{SP} along with the deformation theory of perfect complexes in mixed characteristic as conjectured in \cite[Conjecture 3.1.3]{KLT} and proved recently in \cite{LOfolk}.

\textbf{Claim:} The complex $P$ defines a morphism of stacks $\mu_P: \Y \ra s\mathcal{T}w_{\X_k/k}(0,0,1)^0$ such that at the level of coarse moduli spaces its an open immersion, i.e., $\bar{\mu}: Y \hookrightarrow sTw_{\X_k/k}(0,0,1)^0$ is an open immersion.

This is exactly the content of \cite[Lemma 3.9, 3.11]{Re} (c.f. \cite[Lemma 5.2(ii)]{LO}).s an open immersion. The only thing left to note is that the Mukai vector of $P$ is $(0,0,1)$ since it induces a filtered derived equivalence.

Recall that we can always lift this open immersion uniquely to an open immersion of $Y_A$ into the coarse moduli space of $sTw_{\X_A/A}(0,0,1)^0$, take the open subscheme of $sTw_{\X_A/A}(0,0,1)^0$ supported on $Y$. 

So, now we have the following diagram $$
\xymatrix{
& s\mathcal{T}w_{\X_A/A}(0,0,1)^0 \ar[d]\\
Y_A \ar[r] &sTw_{\X_k/k}(0,0,1)^0}
$$
The vertical arrow on right hand side gives the $\G_m$-gerbe structure to the moduli stack. Now completing the diagram to a cartesian square, gives a $\G_m$-gerbe $\Y_A \ra Y_A$ which is a deformation of $\Y \ra Y$. Restricting the universal complex on $s\mathcal{T}w_{\X_A/A}(0,0,1)^0$ to $\Y_A$ gives us the required kernel, denoted $P_A$. Nakayama lemma implies that $P_A$ also induces a derived equivalence, as in \cite[Theorem 6.1]{LO}.
\end{proof}

\begin{proof}[Proof of Theorem \ref{TwistedpartnersModuli}]
We will the prove this result via lifting to characteristic zero as in \cite[Theorem 6.1, Proposition 8.2]{LO}. 
 
Consider the following composition of twisted derived equivalences
\begin{equation*}
  \Phi_Q:  D^{(-1)}(\Y) \xrightarrow{\Phi_P} D^{(1)}(\X) \xrightarrow{\Phi_{\mathcal{E}}^{-1}} D^{(-1)}(\mathcal{M}_{\X}(v)),  
\end{equation*}
where $\Phi_{\mathcal{E}}^{-1}$ is the inverse of the twisted derived equivalence constructed in Theorem \ref{moduli}. Note that  this composition of twisted derived equivalences is filtered. Indeed, $\Phi_Q((0,0,1)) = (0,0,1)$. Then, the following claim gives us the desired result:

\textbf{Claim:} If $\Phi_Q$ is a filtered derived equivalence, then $\Y \cong \mathcal{M}_{\X}(v)$. 

To unload the notation, we will be working with $Q$ a Fourier-Mukai kernel for a derived equivalence of gerbes $\Y$ and $\X$. Let $H_X$ (resp. $H_Y$) be ample invertible sheaves on $X$ (resp. $Y$), such that $\Phi_Q$ send $H_X$ to $\pm H_Y$. Combining the results of Deligne \cite{LiftingK3} on lifting K3 surfaces and de Jong \cite{deJongGabber} on lifting Brauer classes uniquely (see section \ref{2.2} and \ref{2.4}), we have a lift $(\X_V, H_{X_V})$ of $(\X, H_X)$ over a (possibly) finite extension of the ring of Witt vectors $W(k)$. For all $n \geq 0$, let $V_n := V/(\mathfrak{m}^n)$, where $\mathfrak{m}$ is the maximal ideal of $V$ and let $K := Frac (V)$.  Then Proposition \ref{LOlifting} implies that for each $n$ we have a lifting $\Y_n$ of $\Y$ and a complex $Q_n \in D(\X_n \times \Y_n)$ lifting $Q$. Grothendieck existence theorem for schemes and coherent sheaves and for perfect complexes over algebraic stacks (see, for example \cite{BLim}) gives us a lift $(\Y_V, H_{Y_V})$ over $V$ along with a lift $Q_V$ of $Q$ to $D(\X_V \times \Y_V)$. Using Nakayama lemma, it follows that base change of $Q_V$ to $Q_{K'}$ for any extension $K'$ of $K$, induces a Fourier-Mukai equivalence 
$$
\Phi_{Q_{K'}}: D^b(\Y_{K'}) \xrightarrow{\cong} D^b(\X_{K'}), 
$$
and being filtered implies $\Phi_{Q_{K'}}$ takes $(0,0,1)$ (resp. $(1,0,0)$) to $(0,0,1)$ (resp. $(1,0,0)$). where $(1, 0,0)$ is considered as an element of the extended Mukai lattice. Recall that we have the following commutative diagram of descend to extended Neron-Severi lattice.

\begin{equation*}
   \xymatrixcolsep{10pc}\xymatrix{
    \Phi_P: D^{(-1)}(\mathcal{Y}) \ar[r]^{p_{\X*}(p_{\Y}^*(-)\otimes P)} \ar[d]^{ch_{\Y}(-)} & D^{(1)}(\X) \ar[d]^{ch_{\X}(-)}  \\
\Phi_P^N: \tilde{N}(Y) \ar[r]^{p_{X*}(p_{Y}^*(ch_{\Y}(-)).v_{\X \times \Y}(P)) } &\tilde{N}(X).
}
\end{equation*}
Now choosing an embedding $K \hookrightarrow \C$ gives us a filtered derived equivalence 
$$
\Phi_{Q_{\C}}:D^b(\Y_{\C}) \xrightarrow{\cong} D^b(\X_{\C}),  
$$  
which in turn induces a Hodge isometry 
$$
H^2(Y_{\C}, \Z) \cong H^2(X_{\C}, \Z).
$$
Moreover, by copying the argument in \cite[Step 1 in Proof of Theorem 0.1]{HS2}, we get that $\X_{\C} \cong \Y_{\C}$. Then using the standard argument of spreading out and using the Matsusaka-Mumford result as in \cite[Section 4]{KLT}, we have that the gerbes $\X$ and $\Y$ are isomorphic.  
\end{proof}

\begin{remark}
This can be reformulated in the Caldararu version via \cite[Proposition 2.1.3.3]{L1}. Let $\alpha$ (resp. $\beta$) be the cohomology class of the gerbe $\X \ra X$ (resp. $\mathcal{Y} \ra Y$) in $H^2(X, \G_m) = Br(X)$ (resp. $H^2(Y, \G_m) = Br(Y)$). Then $D^{(-1)}(\X) \cong D^b(X, \alpha^{-1})$ (resp. $D^{(1)}(\mathcal{Y}) \cong D^b(Y, \beta)$) and $D^{(1,1)}(\X \times \mathcal{Y}) \cong D^b(X \times Y, \alpha \boxtimes \beta)$. Thus the above derived equivalence gives us the following derived equivalence:
$$
\Phi_{\mathcal{P}}: D^b(X, \alpha^{-1}) \xrightarrow{\cong} D^b(Y, \beta)
$$
and the content of the above theorem is that $Y$ is isomorphic to the coarse moduli space of $\alpha$-twisted sheaves with Mukai vector $v$, $M_{(X, \alpha)}(v)$ , such that the universal such that the universal sheaf on $M_{(X, \alpha)}(v) \times X$ is a $\alpha \boxtimes \beta$-twisted sheaf.  
\end{remark}

\begin{remark}
Also recall from \cite[2.1.3.10]{Lthesis} that in case $\X \ra X$ is a trivial gerbe or $\mathcal{Y} \ra Y$ is a trivial gerbe or both are trivial gerbes, we have $D^{(1)}(X) \cong D^b(X)$, respectively so in other cases and we recover the partially twisted (or untwisted) statements.    
\end{remark}

\section{Counting Results}

\begin{corollary}[Finiteness of Fourier-Mukai partners]
There  are only finitely many twisted Fourier-Mukai partners of a twisted K3 surface.  
\end{corollary}

\begin{proof}
For K3 surfaces of finite height (including ordinary K3 surfaces) the argument is an extension of the argument in \cite[Section 9]{LO}. Here we the use lifting argument along with the result above that every Fourier-Mukai partner is a moduli space of twisted sheaves. Thus, we lift both of them to characteristic zero and they are bounded by the number of partners of the lift.

For supersingular K3 surfaces, the result is in \cite[Theorem 4.4.6]{Bragg}. 
\end{proof}

\begin{theorem}
The twisted Fourier-Mukai partners of an ordinary K3 surface are in one-to-one correspondence with the twisted Fourier-Mukai partners of the geometric generic fiber of its canonical lift. 
\end{theorem}

\begin{proof}
This is an extension of the result in \cite[Theorem 4.10-4.11]{Sri}.
We again use Lenny's result Theorem \ref{lenny} and Yoshioka's \cite[Theorem 3.19(ii)]{Yos}.
\end{proof}


\section*{Acknowledgements}
Both authors would like to thank Piotr Achinger for helpful mathematical discussions. The first author is grateful to Chennai Mathematical School, India and Tamas Hausel at Institute of Science and Techonology, Austria for supporting her postdoctoral fellowships during which most of the research for the paper was done. This research has received funding from the European Union's Horizon 2020 research and innovation programme under Marie Sklodowska-Curie Grant Agreement No. 754411. The second author was partially supported by the project number 26156 of the Norwegian Research Council and the grant 2019.0493 of the Knut och Alice Wallenberg Stiftelse.

\end{document}